%% file: arxiv.tex
\documentclass[10pt,twoside,reqno,draft]{amsart}
\input{prembu}
\usepackage{amsmath}
\usepackage{amssymb}
\usepackage{amsfonts}
\usepackage{amsthm}
\begin{document}
\setcounter{page}{1}

\leftline{\footnotesize \em {Accepted for Publication in}}
\leftline{\bf\em Fixed Point Theory} \leftline{\footnotesize
http://www.math.ubbcluj.ro/$^{\sim}$nodeacj/sfptcj.html}

\vs*{1.5cm}

\title[Comparable Linear Contractions]{\large Comparable Linear Contractions in Ordered Metric Spaces}
\author[Aftab Alam and Mohammad Imdad]{Aftab Alam$^{*}$ and Mohammad Imdad$^{**}$}
\date{}
\maketitle

\vs*{-0.5cm}

\bc
{\footnotesize
$^{*}$Department of Mathematics, Aligarh Muslim University, Aligarh-202002, India\\
E-mail: aftab\_math.rs@amu.ac.in\\
\medskip
$^{**}$Department of Mathematics, Aligarh Muslim University, Aligarh-202002, India\\
E-mail: mhimdad@amu.ac.in\\
} \ec

\bigskip

{\footnotesize \noindent {\bf Abstract.} In this paper, with a view
to improve the $g$-monotonicity condition, we introduce the notion
of $g$-comparability of a mapping defined on an ordered set and
utilize the same to prove some existence and uniqueness results on
coincidence points for linear contraction without $g$-monotonicity
in ordered metric spaces. Our results extend some classical and well
known results due to Ran and Reurings (Proc. Amer. Math. Soc.
 {\bf 132}(2004), no.5, 1435-1443), Nieto and Rodr\'{\i}guez-L\'{o}pez (Acta Math. Sin. {\bf 23}(2007), no.12,
2205-2212), Turinici (Libertas Math. {\bf 31}(2011), 49-55),
Turinici (Math. Student {\bf 81}(2012), no.1-4, 219-229) and
Dori\'{c} $et\;al.$ (RACSAM {\bf 108}(2014), no.2, 503-510) and
similar others.


}

\bigskip

\section{Introduction}

Throughout this paper, the pair $(X,\preceq),$ stands for a nonempty
set $X$ equipped with a partial order $\preceq$ often called an
ordered set wherein we generally write $x\succeq y$ instead of
$y\preceq x$. Two elements $x$ and $y$ in an ordered set
$(X,\preceq)$ are said to be comparable if either $x\preceq y$ or
$y\preceq x$ and denote it as $x\prec\succ y.$ A subset $E$ of an
ordered set is called totally ordered if $x\prec\succ y$ for all
$x,y\in E.$ In respect of a pair of self-mappings $(f,g)$ defined on
an ordered set $(X,\preceq)$, we say that $f$ is $g$-increasing
(resp. $g$-decreasing) if for any $x,y\in X$, $g(x)\preceq g(y)$
implies $f(x)\preceq f(y) (\rm{resp.}\; f(x)\succeq f(y))$. As per
standard practice, $f$ is called $g$-monotone if $f$ is either
$g$-increasing or $g$-decreasing. Notice that with $g=I$ (the
identity mapping), the notions of $g$-increasing, $g$-decreasing and
$g$-monotone mappings transform into increasing, decreasing and
monotone mappings respectively. Following O'Regan and Petrusel
\cite{PGF2}, the triple $(X,d,\preceq)$ is called ordered metric
space wherein $X$ denotes a nonempty set endowed with a metric $d$
and a partial order $\preceq$. If in addition, $d$ is a complete
metric on $X,$ then we say that $(X,d,\preceq)$ is an ordered
complete
metric space.\\

The relevant detailed discussions on basic topological properties of
ordered sets are available in Milgram \cite{P1,P2}, Eilenberg
\cite{P3}, Wolk \cite{P4,P5} and Monjardet \cite{P6}. Existence of
fixed points for monotone mappings on ordered sets was first
investigated by Tarski \cite{P7} and Bj${\rm\ddot{o}}$rner
\cite{P12} (on complete lattices), Abian and Brown \cite{P19},
DeMarr \cite{P8}, Wong \cite{P9}, Pasini \cite{P20}, Kurepa
\cite{P10}, Amann \cite{P17} and Dugundji and Granas \cite{P21} (on
abstract ordered sets), Ward \cite{P13} (on ordered topological
spaces), DeMarr \cite{P14}(on ordered spaces obtained from complete
metric spaces), Turinici \cite{P15} (on order closed subordered
metrizable uniform spaces) and Turinici \cite{P16} (on quasi-ordered
complete metric spaces). In 2004, unknowingly Ran and Reurings
\cite{PF1} particularized a fixed point theorem proved in Turinici
\cite{P16} in ordered metric spaces for continuous monotone
mappings besides giving some applications to matrix equations.\\

\noindent{\bf Theorem 1.1} (Theorem 2.1, Ran and Reurings
\cite{PF1}). {\it Let $(X,d,\preceq)$ be an ordered metric
space and $f$ a self-mapping on $X$. Suppose that the following conditions hold:\\
\indent\hspace{0.5mm}$(a)$ $(X,d)$ is complete,\\
\indent\hspace{0.5mm}$(b)$ $f$ is monotone,\\
\indent\hspace{0.5mm}$(c)$ $f$ is continuous,\\
\indent\hspace{0.5mm}$(d)$ there exists $x_{0}\in X$ such that
$x_{0}\prec\succ f(x_{0})$,\\
\indent\hspace{0.5mm}$(e)$ there exists $\alpha\in [0,1)$ such that \\
\indent\hspace{2.5cm}$d(fx,fy)\leq\alpha d(x,y)\;\;\forall ~x,y\in X$ with $x\preceq y$,\\
\indent\hspace{0.5mm}$(f)$ every pair of elements of $X$ has a lower bound and an upper bound.\\
Then $f$ has a unique fixed point $\overline{x}$. Moreover, for every $x\in X$, $\lim\limits_{n\to\infty} f^n(x)=\overline{x}.$}\\

 Thereafter,
Nieto and Rodr\'{\i}guez-L\'{o}pez \cite{PF2,PF3} slightly modified
Theorem 1.1 for monotone mappings to relax the continuity
requirement by assuming an additional hypothesis on ordered metric
space besides observing that (owing to assumption $(f)$) the
existence of lower bound(or upper bound) for every pair of elements
of $X$ serves our purpose which is also followed by some
applications of their
results to ordinary differential equations.\\

\noindent{\bf Theorem 1.2} (Theorem 5, Nieto and
Rodr\'{\i}guez-L\'{o}pez \cite{PF3}). {\it Let $(X,d,\preceq)$ be an
ordered metric
space and $f$ a self-mapping on $X$. Suppose that the following conditions hold:\\
\indent\hspace{0.5mm}$(a)$ $(X,d)$ is complete,\\
\indent\hspace{0.5mm}$(b)$ $f$ is monotone,\\
\indent\hspace{0.5mm}$(c)$ either $f$ is continuous or $(X,d,\preceq)$ satisfies the following property:\\
\indent\hspace{0.5cm} if $\{x_n\}$ is a sequence in X such that
$x_n\stackrel{d}{\longrightarrow} x$ whose consecutive terms
are\\\indent\hspace{0.5cm} comparable, then
 there exists a subsequence $\{x_{n_{k}}\}$ of
$\{x_n\}$ such that every\\
\indent\hspace{0.5cm} term is comparable to the limit $x$,\\
\indent\hspace{0.5mm}$(d)$ there exists $x_{0}\in X$ such that
$x_{0}\prec\succ f(x_{0})$,\\
\indent\hspace{0.5mm}$(e)$ there exists $\alpha\in [0,1)$ such that \\
\indent\hspace{2.5cm}$d(fx,fy)\leq\alpha d(x,y)\;\;\forall ~x,y\in X$ with $x\preceq y$,\\
\indent\hspace{0.5mm}$(f)$ every pair of elements of $X$ has a lower
bound or an upper bound.\\
Then $f$ has a unique fixed point.}\\

In subsequent papers (cf.(\cite{PGF1}-\cite{PGF14})) many authors
generalized and refined Theorem 1.2 and proved several fixed point
theorems in ordered metric spaces. In all such results, the
contractivity condition holds on the monotone map for only those
elements which are related by the underlying partial ordering. Thus,
in the context of fixed point theorems for ordered metric spaces,
the usual contraction
condition is weakened but at the expense of monotonicity of the underlying mapping.\\

To relax the monotonicity requirement on underlying mapping, Nieto
and Rodr\'{\i}guez-L\'{o}pez \cite{PF3} replaced this
condition by preservation of comparable elements  and improved Theorem 1.2 as follows:\\

\noindent{\bf Theorem 1.3} (Theorem 7, Nieto and
Rodr\'{\i}guez-L\'{o}pez \cite{PF3}). {\it Let $(X,d,\preceq)$ be an
ordered metric
space and $f$ a self-mapping on $X$. Suppose that the following conditions hold:\\
\indent\hspace{0.5mm}$(a)$ $(X,d)$ is complete,\\
\indent\hspace{0.5mm}$(b)$ for $x,y\in X$ with $x\preceq y\Rightarrow f(x)\preceq f(y)\;{\rm or}\; f(x)\succeq f(y)$ ,\\
\indent\hspace{0.5mm}$(c)$ either $f$ is continuous or $(X,d,\preceq)$ satisfies the following property:\\
\indent\hspace{0.5cm} if $\{x_n\}$ is a sequence in X such that
$x_n\stackrel{d}{\longrightarrow} x$ whose consecutive terms
are\\\indent\hspace{0.5cm} comparable, then
 there exists a subsequence $\{x_{n_{k}}\}$ of
$\{x_n\}$ such that every \\
\indent\hspace{0.5cm} term is comparable to the limit $x$,\\
\indent\hspace{0.5mm}$(d)$ there exists $x_{0}\in X$ such that
$x_{0}\prec\succ f(x_{0})$,\\
\indent\hspace{0.5mm}$(e)$ there exists $\alpha\in [0,1)$ such that \\
\indent\hspace{2.5cm}$d(fx,fy)\leq\alpha d(x,y)\;\;\forall ~x,y\in X$ with $x\preceq y$,\\
\indent\hspace{0.5mm}$(f)$ for every pair $x,y\in X$ there exists
$z\in X$ which is comparable to $x$ and $y$.\\
Then $f$ has a unique fixed point $\overline{x}$. Moreover, for every $x\in X$, $\lim\limits_{n\to\infty} f^n(x)=\overline{x}.$}\\

Here it is noticed that the assumptions $(f)$ of Theorem 1.2 and
$(f)$ of Theorem 1.3 are equivalent (see \cite{PF2}).\\

Turinici \cite{T-RRF,T-NLF} proved similar results besides observing
that these results (hence Theorems 1.1 and 1.2) are particular cases
of Banach Contraction Principle (cf.\cite{B}) and its an important
generalization due to Maia \cite{B6}. Following Turinici
\cite{T-RRF,T-NLF}, given $x,y\in X$, any subset
$\{z_1,z_2,...,z_k\}$ (for $k\geq2$) in $X$ with $z_1=x, z_k =y$ and
$z_i\prec\succ z_{i+1}$ for each $i\in \{1,2,...,{k-1}\}$ is called
a $\prec\succ$-chain between $x$ and $y$. The class of such chains is denoted by C$(x,y,\prec\succ)$.\\

\noindent{\bf Theorem 1.4} (Theorem 2.1, Turinici \cite{T-RRF}).
{\it Let $(X,d,\preceq)$ be an ordered metric
space and $f$ a self-mapping on $X$. Suppose that the following conditions hold:\\
\indent\hspace{0.5mm}$(a)$ $(X,d)$ is complete,\\
\indent\hspace{0.5mm}$(b)$ for $x,y\in X$ with $x\prec\succ y\Rightarrow f(x)\prec\succ f(y)$,\\
\indent\hspace{0.5mm}$(c)$ $f$ is continuous,\\
\indent\hspace{0.5mm}$(d)$ there exists $x_{0}\in X$ such that
$x_{0}\prec\succ f(x_{0})$,\\
\indent\hspace{0.5mm}$(e)$ there exists $\alpha\in [0,1)$ such that \\
\indent\hspace{2.5cm}$d(fx,fy)\leq\alpha d(x,y)\;\;\forall ~x,y\in X$ with $x\preceq y$,\\
\indent\hspace{0.5mm}$(f)$ C$(x,y,\prec\succ)$ is nonempty for each
$x,y\in X$.\\
Then $f$ has a unique fixed point $z$. Moreover for each $x\in X$, the sequence $\{f^nx\}$ is convergent and $\lim\limits_{n\to\infty} f^n(x)=z.$}\\

\noindent{\bf Theorem 1.5} (Theorem 2.1, Turinici \cite{T-NLF}).
{\it Let $(X,d,\preceq)$ be an ordered metric
space and $f$ a self-mapping on $X$. Suppose that the following conditions hold:\\
\indent\hspace{0.5mm}$(a)$ $(X,d)$ is complete,\\
\indent\hspace{0.5mm}$(b)$ for $x,y\in X$ with $x\prec\succ y\Rightarrow f(x)\prec\succ f(y)$,\\
\indent\hspace{0.5mm}$(c)$  $(X,d,\preceq)$ satisfies the following property:\\
\indent\hspace{0.5cm} if $\{x_n\}$ is a sequence in X such that
$x_n\stackrel{d}{\longrightarrow} x$ whose consecutive terms
are\\\indent\hspace{0.5cm} comparable, then
 there exists a subsequence $\{x_{n_{k}}\}$ of
$\{x_n\}$ such that every\\
\indent\hspace{0.5cm} term is comparable to the limit $x$,\\
\indent\hspace{0.5mm}$(d)$ there exists $x_{0}\in X$ such that
$x_{0}\prec\succ f(x_{0})$,\\
\indent\hspace{0.5mm}$(e)$ there exists $\alpha\in [0,1)$ such that \\
\indent\hspace{2.5cm}$d(fx,fy)\leq\alpha d(x,y)\;\;\forall ~x,y\in X$ with $x\preceq y$,\\
\indent\hspace{0.5mm}$(f)$ C$(x,y,\prec\succ)$ is nonempty for each
$x,y\in X$.\\
Then $f$ has a unique fixed point $z$. Moreover for each $x\in X$, the sequence $\{f^nx\}$ is convergent and $\lim\limits_{n\to\infty} f^n(x)=z.$}\\

Notice that the assumptions $(b)$ of Theorems 1.4 and 1.5 are
equivalents to assumption $(b)$ of Theorems 1.3. But assumptions
$(f)$ of Theorems 1.4 and 1.5 are relatively weaker than assumption
$(f)$ of
Theorems 1.3 (see details in \cite{T-RRF,T-NLF}).\\

Very recently Dori\'{c} $et\;al.$
\cite{PF-C0} proved the following result:\\

\noindent{\bf Theorem 1.6} (Corollary 2.7, Dori\'{c} $et\;al.$
\cite{PF-C0}). {\it Let $(X,d,\preceq)$ be an ordered metric
space and $f$ a self-mapping on $X$. Suppose that the following conditions hold:\\
\indent\hspace{0.5mm}$(a)$ $(X,d)$ is complete,\\
\indent\hspace{0.5mm}$(b)$ for $x,y\in X$ with $x\prec\succ y\Rightarrow f(x)\prec\succ f(y)$,\\
\indent\hspace{0.5mm}$(c)$  $f$ is continuous or $(X,d,\preceq)$ satisfies the following property:\\
\indent\hspace{0.5cm} if $x_n\stackrel{d}{\longrightarrow} x$ in $X$
then $x_n\prec\succ x$ for $n$ sufficiently large,\\
\indent\hspace{0.5mm}$(d)$ there exists $x_{0}\in X$ such that
$x_{0}\prec\succ f(x_{0})$,\\
\indent\hspace{0.5mm}$(e)$ there exists $\alpha\in [0,1)$ such that \\
\indent\hspace{2.5cm}$d(fx,fy)\leq\alpha d(x,y)\;\;\forall ~x,y\in X$ with $x\prec\succ y$.\\
Then $f$ has a fixed point in $X$. Moreover if\\
\indent\hspace{0.5mm}$(f)$ for each $x,y\in X\;\exists~z\in X$ such
that $x\prec\succ z$ and $y\prec\succ z$.\\ Then the fixed point
$u$ is unique and for each $x\in X$ sequence $\{f^nx\}$ converges to $u$.}\\

Although the property on $(X,d,\preceq)$ in assumption $(c)$ of
Theorem 1.3 is relatively weaker than the assumption $(c)$ of
Theorem 1.6. But in Theorem 1.6, authors observed that the
uniqueness of fixed point is not necessary.\\

As reflected in Theorems 1.3-1.6, with a view to coin a relatively
weaker alternate condition to avoid the use of the monotonicity
requirement, the respective authors used a common property on the
involved mapping (see assumption $(b)$ in Theorems 1.3-1.6). In this
paper, we generalize this idea to a pair of mappings and utilize the
same to prove some coincidence point theorems for a pair of
self-mappings $f$ and $g$ defined on an ordered metric space $X$
satisfying linear contractivity condition in two different
directions namely: in case $X$ is complete or alternately $X$ has a
complete subspace $Y$ such that $f(X)\subseteq Y\subseteq g(X)$,
while the whole space $X$ may or may not be complete. As a
consequence of our results, we also derive a corresponding fixed
point theorem, which extends and improves all earlier mentioned
results ($i.e.$ Theorems 1.1-1.6) besides
furnishing an illustrative example.\\

\bigskip

\section{Preliminaries}
In this section, we summarize some basic definitions and auxiliary
results. Throughout this paper, $\mathbb{N}$ stands for the set of
natural numbers, while
$\mathbb{N}_0$ for the set of whole numbers ($i.e.$ $\mathbb{N}_0=\mathbb{N}\cup \{0\}$).\\

\noindent{\bf Definition 2.1 \cite{CP1,CP2}}. Let $X$ be a nonempty
set and $(f,g)$ a pair of self-mappings on $X$. Then
\begin{enumerate}
\item[{(i)}] an element $x\in X$ is called a
coincidence point of $f$ and $g$ if
$$g(x)=f(x),$$
\item[{(ii)}] if $x\in X$ is a coincidence point of $f$ and $g$
and $\overline{x}\in X$ such that $\overline{x}=g(x)=f(x)$, then
$\overline{x}$ is called a point of coincidence of $f$ and $g$,
\item[{(iii)}] if $x\in X$ is a coincidence point of $f$ and $g$
such that $x=g(x)=f(x)$, then $x$ is called a common fixed point of
$f$ and $g$,
\item[{(iv)}] the pair $(f,g)$ is said to be commuting if
$$g(fx)=f(gx)\hspace{0.5cm}\forall~x\in X\;{\rm and}$$
\item[{(v)}] the pair
$(f,g)$ is said to be weakly compatible (or partially commuting or
coincidentally commuting) if $f$ and $g$ commute at their
coincidence points, $i.e.$,
$$g(fx)=f(gx)\hspace{0.5cm}{\rm whenever}\;\; g(x)=f(x).$$
\end{enumerate}

\noindent{\bf Definition 2.2 \cite{CP3,CP4}}. Let $(X,d)$ be a
metric space and $(f,g)$ a pair of self-mappings on $X$. Then
\begin{enumerate}
\item[{(i)}] the pair
$(f,g)$ is said to be weakly commuting if
$$d(gfx,fgx)\leq d(gx,fx)\hspace{0.3cm}\forall~x\in X\;{\rm and}$$
\item[{(ii)}] the pair $(f,g)$ is said to be compatible if
$$\lim\limits_{n\to \infty}d(gfx_n,fgx_n)=0$$
whenever $\{x_n\}$ is a sequence in $X$ such that
$$\lim\limits_{n\to \infty}g(x_n)=\lim\limits_{n\to \infty}f(x_n).$$
\end{enumerate}

Clearly in a metric space, commutativity $\Rightarrow$ weak
commutativity $\Rightarrow$ compatibility $\Rightarrow$ weak
compatibility but reverse implications are not true in
general (for details see \cite{CP2}-\cite{CP4}).\\

\noindent{\bf Definition 2.3 \cite{g-C}}. Let $(X,d)$ be a metric
space, $(f,g)$ a pair of self-mappings on $X$ and $x\in X$. We say
that $f$ is $g$-continuous at $x$ if for all $\{x_n\}\subset X$,
$$g(x_n)\stackrel{d}{\longrightarrow} g(x)\Rightarrow f(x_n)\stackrel{d}{\longrightarrow} f(x).$$
Moreover, $f$ is called $g$-continuous if it is $g$-continuous at
each
point of $X$.\\
Notice that with $g=I$ (the identity mapping on $X$)
Definition 2.3 reduces to the definition of continuity.\\

\noindent{\bf Definition 2.4 \cite{PGF14}}. Let $(X,\preceq)$ be an
ordered set and $\{x_n\}\subset X.$ Then
\begin{enumerate}
\item[{(i)}] the sequence $\{x_n\}$ is said to be termwise bounded if there is an element $z\in X$ such that each term
of $\{x_n\}$ is comparable with $z,$ $i.e.$,
$$x_n\prec\succ z\;\;\;\;\;\;\;\forall~ n\in \mathbb{N}_0$$
so that $z$ is a c-bound of $\{x_n\}$ and
\item[{(ii)}] the sequence $\{x_n\}$ is said to be termwise monotone if consecutive terms of $\{x_n\}$ are comparable,
$i.e.$,
$$x_n\prec\succ x_{n+1}\;\;\forall~ n\in \mathbb{N}_0.$$
\end{enumerate}

Clearly all bounded above as well as bounded below sequences are
termwise bounded and all monotone sequences are termwise monotone.\\

Let $(X,d,\preceq)$ be an ordered metric space and $\{x_n\}\subset
X$. If $\{x_n\}$ is termwise monotone and
$x_n\stackrel{d}{\longrightarrow} x$, then we denote it
symbolically by $x_n\updownarrow x.$\\

\noindent{\bf Definition 2.5 \cite{PGF14}}. Let $(X,d,\preceq)$ be
an ordered metric space. We say that $(X,d,\preceq)$ has {\it
TCC}\;(termwise monotone-convergence-c-bound) property if every
termwise monotone convergent sequence $\{x_n\}$ in $X$ has a
subsequence, which is termwise bounded by the limit of $\{x_n\}$ (as
a c-bound), $i.e.,$
$$x_n\updownarrow x \Rightarrow
\exists {\rm \;a\; subsequence\;} \{x_{n_k}\}{\rm \;of\;} \{x_n\}
\;{\rm with}\;\;x_{n_k}\prec\succ x\;\; \forall~ k\in
\mathbb{N}_0.$$

Notice that above definition is formulated using a property utilized
in Nieto and Rodr\'{\i}guez-L\'{o}pez \cite{PF3} (see assumption
$(c)$ in Therem 1.2).\\

\noindent{\bf Definition 2.6 \cite{PGF14}}. Let $(X,d,\preceq)$ be
an ordered metric space and $g$ a self-mapping on $X.$ We say that
$(X,d,\preceq)$ has {\it g-TCC} property if every termwise monotone
convergent sequence $\{x_n\}$ in $X$ has a subsequence, whose
$g$-image is termwise bounded by $g$-image of limit of $\{x_n\}$ (as
a c-bound),
 $i.e.,$
$$x_n\updownarrow x \Rightarrow \exists
{\rm \;a\; subsequence\;} \{x_{n_k}\}{\rm \;of\;} \{x_n\} \;{\rm
with}\;\;g(x_{n_k})\prec\succ g(x)\;\; \forall~ k\in \mathbb{N}_0.$$
Notice that under the restriction $g=I,$ the identity
mapping on $X,$ Definition 2.6 reduces to Definition 2.5.\\

Very recently, Alam and Imdad \cite{PGF14} generalized Theorems 1.1
and 1.2 for a pair of mappings and proved the following coincidence
point theorem:\\

\noindent{\bf Theorem 2.7 \cite{PGF14}.} {\it Let $(X,d,\preceq)$ be
an ordered
metric space and $f$ and $g$ two self-mappings on $X$. Suppose that the following conditions hold:\\
\indent\hspace{0.5mm}$(a)$ $f(X)\subseteq g(X)$,\\
\indent\hspace{0.5mm}$(b)$ $f$ is $g$-monotone,\\
\indent\hspace{0.5mm}$(c)$ there exists $x_{0}\in X$ such that $g(x_{0})\prec\succ f(x_{0})$,\\
\indent\hspace{0.5mm}$(d)$ there exists $\alpha\in [0,1)$ such that \\
\indent\hspace{2.5cm}$d(fx,fy)\leq\alpha d(gx,gy)\;\;\forall~ x,y\in X$ with $g(x)\prec\succ g(y)$,\\
\indent\hspace{0.5mm}$(e)$ $(e1)$ $(X,d)$ is complete,\\
\indent\hspace{7mm}$(e2)$ $(f,g)$ is compatible pair,\\
\indent\hspace{7mm}$(e3)$ $g$ is continuous,\\
\indent\hspace{7mm}$(e4)$ either $f$ is continuous or $(X,d,\preceq)$ has {\it g-TCC} property,\\
\indent\hspace{1mm}or alternately\\
\indent\hspace{0.5mm}$(e^\prime)$ $(e^{\prime}1)$ either $(fX,d)$ or $(gX,d)$ is complete,\\
\indent\hspace{7mm}$(e^{\prime}2)$ either $f$ is $g$-continuous or $f$ and $g$ are continuous or $(gX,d,\preceq)$ has\\
\indent\hspace{1.5cm}{\it TCC} property.\\
Then $f$ and $g$ have a coincidence point.}\\

We need the following known results in the proof of our main results. \\

\noindent{\bf Lemma 2.8 \cite{L2}.} {\it Let $X$ be a nonempty set
and $g$ a self-mapping on $X$. Then there exists a subset
$E\subseteq X$ such that $g(E)=g(X)$ and
$g: E \rightarrow X$ is one-one.}\\

\noindent{\bf Lemma 2.9 \cite{PGF13}.} {\it Let $X$ be a nonempty
set and $f$ and $g$ two self-mappings on $X$. If the pair $(f,g)$ is
weakly compatible, then every point of coincidence of $f$ and $g$ is
also a coincidence point of $f$ and $g$.}\\
\bigskip
\section{Results on Coincidence Points}
\noindent Firstly, we name a property utilized in Theorems 1.3-1.6
and term the same as comparable mapping.\\

\noindent{\bf Definition 3.1} (see \cite{PF3,T-RRF,T-NLF,PF-C0}).
Let $(X,\preceq)$ be an ordered set and $f$ a self-mapping on $X$.
We say that $f$ is comparable (or weakly monotone or
$\prec\succ$-preserving) if $f$ maps comparable elements to
comparable elements, $i.e.$, for any $x,y\in X$
$$x\prec\succ y\Rightarrow f(x)\prec\succ f(y).$$
It is clear that every monotone mapping is comparable, but not
conversely. To substantiate this view point, consider the set
$X=[-\frac{1}{3},\frac{1}{3}]$ under the natural ordering of real
numbers. Define $f:X\rightarrow X$ by $f(x)=x^2$, then $f$ is
comparable
but not monotone.\\

We extend the idea embodied in Definition 3.1 to a pair of mappings
to introduce
the notion of $g$-comparability:\\

\noindent{\bf Definition 3.2}. Let $(X,\preceq)$ be an ordered set
and $f$ and $g$ two self-mappings on $X$. We say that $f$ is
$g$-comparable (or weakly $g$-monotone or
$(g,\prec\succ)$-preserving) if for any $x,y\in X$
$$g(x)\prec\succ g(y)\Rightarrow f(x)\prec\succ f(y).$$
Notice that on setting $g=I$, the identity mapping on $X$,
Definition 3.2 reduces to Definition 3.1.\\

Now, we are equipped to prove our main result on coincidence points in ordered complete metric spaces which runs as follows:\\

\noindent{\bf Theorem 3.3.} {\it Let $(X,d,\preceq)$ be an ordered
complete metric space and $f$ and $g$ two self-mappings on $X$. Suppose that the following conditions hold:\\
\indent\hspace{0.5mm}(i) $f(X)\subseteq g(X)$,\\
\indent\hspace{0.5mm}(ii) $f$ is $g$-comparable,\\
\indent\hspace{0.5mm}(iii) $(f,g)$ is compatible pair,\\
\indent\hspace{0.5mm}(iv) $g$ is continuous,\\
\indent\hspace{0.5mm}(v) either $f$ is continuous or $(X,d,\preceq)$ has {\it g-TCC} property,\\
\indent\hspace{0.5mm}(vi) there exists $x_{0}\in X$ such that $g(x_{0})\prec\succ f(x_{0})$,\\
\indent\hspace{0.5mm}(vii) there exists $\alpha\in [0,1)$ such that \\
\indent\hspace{2.5cm}$d(fx,fy)\leq\alpha d(gx,gy)\;\;\forall~ x,y\in
X$ with $g(x)\prec\succ g(y)$.\\
Then $f$ and $g$ have a coincidence point.}\\
{\bf Proof.} In view of assumption (vi) if $g(x_0)=f(x_0)$, then we
are through. Otherwise, if $g(x_0)\neq f(x_0)$, then owing to
assumption (i) ($i.e.$ $f(X)\subseteq g(X)$), we can choose $x_1\in
X$ such that $g(x_1)=f(x_0)$. Again from $f(X)\subseteq g(X)$, we
can choose $x_2\in X$ such that $g(x_2)=f(x_1)$. Continuing this
process inductively, we define a sequence $\{x_n\}\subset X$ of
joint iterates such that
$$g(x_{n+1})=f(x_n)\;\;\forall ~n \in \mathbb{N}_0.\eqno (1)$$
Now, we assert that $\{gx_n\}$ is a termwise monotone sequence,
$i.e.,$
$$g(x_n)\prec\succ g(x_{n+1})\;\;\forall~ n \in \mathbb{N}_0.\eqno (2)$$
We prove this fact by mathematical induction. On using assumption
(vi) and equation (1) with $n=0$, we have
$$g(x_{0})\prec\succ f(x_{0})=g(x_1)$$
Thus, (2) holds for $n=0.$ Suppose that (2) holds for  $n=r>0$,
$i.e.,$
$$g(x_r)\prec\succ g(x_{r+1})\eqno (3)$$
then we have to show that (2) holds for $n=r+1$. To accomplish this,
we use (1), (3) and assumption (ii) so that
$$g(x_{r+1})=f(x_r)\prec\succ f(x_{r+1})=g(x_{r+2}).$$
Thus, by induction, (2) holds for all $n \in \mathbb{N}_0$.\\
If $g(x_{n_0})=g(x_{n_0+1})$ for some $n_0\in \mathbb{N}$, then
using (1), we have $g(x_{n_0})=f(x_{n_0})$, $i.e.,$ $x_{n_0}$ is a
coincidence point of $f$ and $g$ and hence we are done. On the other
hand, if $g(x_n)\neq g(x_{n+1})$ for each $n\in \mathbb{N}_0$, then
$d(gx_n,gx_{n+1})\neq 0$ for each $n\in \mathbb{N}_0$. On using (1),
(2) and assumption (vii), we obtain
$$d(gx_{n},gx_{n+1})=d(fx_{n-1},fx_{n})\leq \alpha d(gx_{n-1},gx_{n})\;\;\;\forall~n\in \mathbb{N}.$$
By induction, we have
$$d(gx_{n},gx_{n+1})\leq \alpha d(gx_{n-1},gx_{n})\leq \alpha^2 d(gx_{n-2},gx_{n-1})\leq\cdots\leq \alpha^n d(gx_{0},gx_{1})\;\;\forall~n\in\mathbb{N}$$
so that
$$d(gx_{n},gx_{n+1})\leq \alpha^n d(gx_{0},gx_{1})\;\;\;\forall~n\in \mathbb{N}.\eqno(4)$$
For $n<m$, using (4), we obtain
\begin{eqnarray*}
d(gx_{n},gx_{m})&\leq&
d(gx_{n},gx_{n+1})+d(gx_{n+1},gx_{n+2})+\cdots+d(gx_{m-1},gx_{m})\\
&\leq& (\alpha^n+\alpha^{n+1}+\cdots+\alpha^{m-1})d(gx_{0},gx_{1})\\
&=& \frac{\alpha^n-\alpha^m}{1-\alpha}d(gx_{0},gx_{1})\\
&\leq& \frac{\alpha^n}{1-\alpha}d(gx_{0},gx_{1})\\
&\rightarrow& 0\;{\rm as}\;m,n\rightarrow \infty.
\end{eqnarray*}
Therefore $\{gx_n\}$ is a Cauchy sequence. As $X$ is complete, there
exists $z\in X$ such that
$$\lim\limits_{n\to\infty} g(x_n)=z.\eqno(5) $$
On using (1) and (5), we obtain
$$\lim\limits_{n\to\infty} f(x_n)=\lim\limits_{n\to\infty} g(x_{n+1})=z.\eqno(6) $$
On using continuity of $g$ in (5) and (6), we get
$$\lim\limits_{n\to\infty} g(gx_n)=g(\lim\limits_{n\to\infty}gx_n)=g(z).\eqno(7) $$
$$\lim\limits_{n\to\infty} g(fx_n)=g(\lim\limits_{n\to\infty}fx_n)=g(z).\eqno(8) $$
As $\lim\limits_{n\to\infty} f(x_n)=\lim\limits_{n\to\infty}
g(x_n)=z$ (due to (5) and (6)), on using compatibility of $f$ and
$g$, we obtain
$$\lim\limits_{n\to\infty}d(gfx_n,fgx_n)=0.\eqno(9)$$
Now, we show that $z$ is a coincidence point of $f$ and $g$. To
accomplish this, we use assumption (v). Suppose that $f$ is
continuous. On using (5) and continuity of $f$, we obtain
$$\lim\limits_{n\to\infty} f(gx_n)=f(\lim\limits_{n\to\infty} gx_n)=f(z).\eqno(10) $$
On using (8), (9), (10) and continuity of $d$, we obtain
\begin{eqnarray*}
d(gz,fz)&=&d(\lim\limits_{n\to\infty}gfx_{n},\lim\limits_{n\to\infty}fgx_{n})\\
&=&\lim\limits_{n\to\infty}d(gfx_n,fgx_n)\\
&=&0
\end{eqnarray*}
so that $$ g(z)=f(z).$$ Thus $z$ is a coincidence point of $f$
and $g$ and hence we are through.\\
Alternately, suppose that $(X,d,\preceq)$ has $g$-{\it TCC}
property. As $g(x_n)\updownarrow z$ (due to (2) and (5)), $\exists~$
a subsequence $\{y_{n_k}\}$ of $\{gx_n\}$ such that
$$g(y_{n_k})\prec\succ g(z)\;\forall~ k\in \mathbb{N}_0.$$
Now $\{gx_n\}\subset g(X)$ and $\{y_{n_k}\}\subset \{gx_n\}, \exists
\{x_{n_k}\}\subset X$ such that $y_{n_k}=g(x_{n_k})$. Hence, we have
$$g(gx_{n_k})\prec\succ g(z)\;\forall~ k\in
\mathbb{N}_0.\eqno(11)$$ Since $g(x_{n_k})\rightarrow z,$ so
equations (5)-(10) hold for also $\{x_{n_k}\}$ instead of $\{x_n\}$.
On using (11) and assumption (vii), we obtain
$$d(fgx_{n_k},fz)\leq\alpha d(ggx_{n_k},gz)\;\;\forall~ k\in \mathbb{N}_0.\eqno(12)$$
On using triangular inequality, (7), (8), (9) and (12), we get
\begin{eqnarray*}
 \nonumber d(gz,fz)&\leq& d(gz,gfx_{n_k})+d(gfx_{n_k},fgx_{n_k})+d(fgx_{n_k},fz)
\\&\leq&d(gz,gfx_{n_k})+d(gfx_{n_k},fgx_{n_k})+\alpha d(ggx_{n_k},gz)\\
&\rightarrow& 0\;{\rm as}\; k\rightarrow \infty
\end{eqnarray*}
so that $$ g(z)=f(z).$$
Thus $z$ is a coincidence point of $f$ and $g$ and hence this concludes the proof.\\

As commutativity $\Rightarrow$ weak commutativity $\Rightarrow$
compatibility for a pair of mappings, therefore the following
consequence of Theorem 3.3 trivially holds.\\

\noindent{\bf Corollary 3.4.} {\it Theorem 3.3 remains true if we
replace condition (iii) by one of the following conditions besides retaining the rest of the hypotheses:\\
\indent\hspace{0.5mm}(iii)$^\prime$ $(f,g)$ is commuting pair,\\
\indent\hspace{0.5mm}(iii)$^{\prime\prime}$ $(f,g)$ is weakly
commuting pair.}\\

Our next result is analogous to Theorem 3.3 whenever $X$ is not
necessarily complete. Instead, we require at least, one
of its subspaces to be complete.\\

\noindent{\bf Theorem 3.5.} {\it Let $(X,d,\preceq)$ be an ordered
metric space, $Y\subseteq X$ and $f$ and $g$ two self-mappings on $X$. Suppose that the following conditions hold:\\
\indent\hspace{0.5mm}(i) $f(X)\subseteq Y \subseteq g(X)$,\\
\indent\hspace{0.5mm}(ii) $f$ is $g$-comparable,\\
\indent\hspace{0.5mm}(iii) $(Y,d)$ is complete,\\
\indent\hspace{0.5mm}(iv) either $f$ is $g$-continuous or $f$ and $g$ are continuous or $(Y,d,\preceq)$ has {\it TCC}\\
\indent\hspace{7mm}property,\\
\indent\hspace{0.5mm}(v) there exists $x_{0}\in X$ such that $g(x_{0})\prec\succ f(x_{0})$,\\
\indent\hspace{0.5mm}(vi) there exists $\alpha\in [0,1)$ such that \\
\indent\hspace{2.5cm}$d(fx,fy)\leq\alpha d(gx,gy)\;\;\forall~ x,y\in X$ with $g(x)\prec\succ g(y)$.\\
Then $f$ and $g$ have a coincidence point.}\\
{\bf Proof.} In view of assumption (v) if $g(x_0)=f(x_0)$, then
$x_0$ is a coincidence point of $f$ and $g$ and hence proof is
completed. Otherwise, if $g(x_0)\neq f(x_0)$, then we have
$g(x_0)\prec f(x_0)$. As $g$ is a self-mapping on $X$, by using
Lemma 2.8, there exists a subset $E\subseteq X$ such that
$g(E)=g(X)$ and $g: E \rightarrow X$ is one-one. Hence assumption
(i) implies that $f(X)\subseteq g(E)$ so that we can choose $e_1\in
E$ such that $g(e_1)=f(x_0)$. Again, we can choose $e_2\in E$ such
that $g(e_2)=f(e_1)$. Now proceeding in the same way, we can
inductively construct a sequence $\{e_n\}\subset E$ such that
$$g(e_{n+1})=f(e_n)\;\;\forall ~n \in \mathbb{N}.\eqno (13)$$ Following the proof of Theorem
3.3, we can show that the sequence $\{ge_n\}$ (and hence $\{fe_n\}$
also) is termwise monotone and also Cauchy. Owing to $f(X)\subseteq
Y$, $\{fe_n\}$ is a Cauchy sequence in $Y$. As $Y$ is complete,
there exists $y\in Y$ such that $\lim\limits_{n\to\infty} f(e_n)=y$.
Now, by assumption $Y\subseteq g(X)=g(E)$, we can find $e\in E$ such
that $y=g(e)$. Hence, on using (13), we have
$$\lim\limits_{n\to\infty} g(e_n)=\lim\limits_{n\to\infty} f(e_n)=g(e).\eqno(14)$$
Now, we show that $e$ is a coincidence point of $f$ and $g$. In view
of assumption (iv), firstly, suppose that $f$ is $g$-continuous,
then using (14), we get
$$\lim\limits_{n\to\infty} f(e_n)=f(e).\eqno(15)$$
On using (14), (15) and uniqueness of limit, we get
$$g(e)=f(e),$$
$i.e.$, $e$ is a coincidence point of $f$ and $g$ and hence we are through.\\
Secondly, suppose that $f$ and $g$ both are continuous. Define $T:
g(E) \rightarrow g(E)$ by
$$T(ga)=f(a)\;\;\forall\; g(a)\in g(E)\; {\rm where}\; a\in E.\eqno(16)$$
As $g: E \rightarrow X$ is one-one and $f(X)\subseteq g(E)$, $T$ is
well defined. Again since $f$ and $g$ are continuous, it follows
that $T$ is continuous. On using (14), (16) and continuity of $T$,
we get
$$f(e)=T(ge)=T(\lim\limits_{n\to\infty} ge_n)=\lim\limits_{n\to\infty} T(ge_n)=\lim\limits_{n\to\infty} f(e_n)=g(e).$$
Thus $e$ is a coincidence point of $f$ and $g$ and hence we are done.\\
Finally, suppose that $(Y,d,\preceq)$ has {\it TCC} property. Using
(13) and assumption (i), $\{ge_n\}$ is termwise monotone in $Y$ and
using (14), $g(e_n)\stackrel{d}{\longrightarrow} g(e)$, which yield
that $g(e_n)\updownarrow g(e)$. Hence, by {\it TCC} property of $Y$,
$\exists~$ a subsequence $\{ge_{n_k}\}$ of $\{ge_n\}$ such that
$$g(e_{n_k})\prec\succ g(e)\;\forall~ k\in \mathbb{N}.\eqno(17)$$
On using (17) and assumption (vi), we obtain
$$d(fe_{n_k},fe)\leq \alpha d(ge_{n_k},ge)\;\;\forall~ k\in \mathbb{N}.\eqno(18)$$
On using (14), (18) and continuity of $d$, we get
\begin{eqnarray*}
 \nonumber d(ge,fe)&=&d(\lim\limits_{k\to\infty} fe_{n_k},fe)
\\&=&\lim\limits_{k\to\infty}d(fe_{n_k},fe)
\\&\leq&\alpha\lim\limits_{k\to\infty}d(ge_{n_k},ge)\\
&=& 0
\end{eqnarray*}
so that $$g(e)=f(e).$$
Hence $e$ is a coincidence point of $f$ and $g$. This completes the proof.\\

Now, we present a consequence of Theorem 3.5.\\

\noindent{\bf Corollary 3.6.} {\it Theorem 3.5 remains true if we
replace (iii) by one of the following
conditions (iii)$^\prime$ and (iii)$^{\prime\prime}$ besides retaining the rest of the hypotheses:\\
\indent\hspace{0.5mm}(iii)$^\prime$ $(X,d)$ is complete and
one of $f$ and $g$ is onto,\\
\indent\hspace{1mm}(iii)$^{\prime\prime}$ $(X,d)$ is
complete and $Y$ (where $f(X)\subseteq Y \subseteq g(X)$) is a closed subspace.}\\
{\bf Proof.} If (iii)$^\prime$ holds, we get either $f(X)=X$ or
$g(X)=X$ so that either $f(X)$ or $g(X)$ is complete, which implies
that (iii) holds and hence Theorem 3.5 is
applicable.\\
If (iii)$^{\prime\prime}$ holds, then using the fact that closed
subset of a complete metric space is complete, $Y$ is complete which
implies that (iii) holds and hence Theorem 3.5 is
applicable.\\

On combining Theorems 3.3 and 3.5, we obtain the following result:\\

\noindent{\bf Theorem 3.7.} {\it Let $(X,d,\preceq)$ be an ordered
metric space and $f$ and $g$ two self-mappings on $X$. Suppose that the following conditions hold:\\
\indent\hspace{0.5mm}$(a)$ $f(X)\subseteq g(X)$,\\
\indent\hspace{0.5mm}$(b)$ $f$ is $g$-comparable,\\
\indent\hspace{0.5mm}$(c)$ there exists $x_{0}\in X$ such that $g(x_{0})\prec\succ f(x_{0})$,\\
\indent\hspace{0.5mm}$(d)$ there exists $\alpha\in [0,1)$ such that \\
\indent\hspace{2.5cm}$d(fx,fy)\leq\alpha d(gx,gy)\;\;\forall~ x,y\in X$ with $g(x)\prec\succ g(y)$,\\
\indent\hspace{0.5mm}$(e)$ $(e1)$ $(X,d)$ is complete,\\
\indent\hspace{7mm}$(e2)$ $(f,g)$ is compatible pair,\\
\indent\hspace{7mm}$(e3)$ $g$ is continuous,\\
\indent\hspace{7mm}$(e4)$ either $f$ is continuous or $(X,d,\preceq)$ has {\it g-TCC} property,\\
\indent\hspace{2mm}or alternately\\
\indent\hspace{0.5mm}$(e^\prime)$ $(e^{\prime}1)$ there exists a
subset $Y$ of $X$ such that $f(X)\subseteq Y \subseteq g(X)$ and
$(Y,d)$ is\\
\indent\hspace{1.5cm}complete,\\
\indent\hspace{7mm}$(e^{\prime}2)$ either $f$ is $g$-continuous or $f$ and $g$ are continuous or $(Y,d,\preceq)$ has\\
\indent\hspace{1.5cm}{\it TCC} property.\\
Then $f$ and $g$ have a coincidence point.}\\

Notice that Theorem 3.7 improves Theorem 2.7 and hence in Theorem
2.7 the $g$-monotonicity can be alternately replaced by
$g$-comparability, which is relatively weaker.\\
\bigskip
\section{Uniqueness Results}
Recall that in order to obtain the uniqueness of fixed point in
ordered metric spaces, several authors used the following
alternative
conditions.\\

(I) $(X,\preceq)$ is totally ordered.\\

The preceding condition is more natural, as under this condition,
results of Ran and Reurings \cite{PF1} and Nieto and
Rodr\'{\i}guez-L\'{o}pez \cite{PF2,PF3} follow directly from Banach
contraction principle \cite{B} but this condition is very restrictive.\\

Ran and Reurings \cite{PF1} used the following condition to obtain
the
uniqueness of fixed point in their result (see Theorem 1.1).\\

(II) every pair of elements of $X$ has a lower bound and an upper
bound.\\

Later, Nieto and Rodr\'{\i}guez-L\'{o}pez \cite{PF2,PF3} (see
Theorem 1.2) modified condition (II) by
assuming relatively weaker condition as follows:\\

(III) every pair of elements of $X$ has a lower bound or an upper
bound,\\

which is equivalent (proved in \cite{PF2}) to the following:\\

\indent\hspace{5mm} for each pair $x,y\in X,\;\exists~z\in X$ such
that $x\prec\succ z$ and
$y\prec\succ z$.\\

On the lines of Jleli $et\;al.$ \cite{PGF10}, $(X,\preceq)$ is
called
directed if it satisfies condition (III) (see Definition 2.4 \cite{PGF10}).\\

Turinici \cite{T-RRF,T-NLF} used the following condition (see Theorems 1.4 and 1.5):\\

(IV) C$(x,y,\prec\succ)$ is nonempty, for each
$x,y\in X$.\\

Clearly, (I)$\Rightarrow$(II)$\Rightarrow$(III)$\Rightarrow$(IV)
$i.e.$ among these four conditions (IV) is the weakest one.\\

Inspired by Jleli $et\;al.$ \cite{PGF10}, we extend condition (III) to a pair of mappings which runs as follows:\\

\noindent{\bf Definition 4.1.} Let $(X,\preceq)$ be an ordered set
and $(f,g)$ a pair of self-mappings on $X$. We say that
$(X,\preceq)$ is $(f,g)$-directed if for each pair $x,y\in X$,
$\exists~z\in X$ such that $f(x)\prec\succ g(z)$ and $f(y)\prec\succ
g(z)$.\\
In cases $f=I$ and $f=g=I$ (where $I$ denotes identity mapping on
$X$), $(X,\preceq)$ is called $g$-directed and directed
respectively.\\

Inspired by Turinici \cite{T-RRF}, we limit condition (IV) to an arbitrary subset rather than the whole ordered set which runs as follows:\\

\noindent{\bf Definition 4.2.} Let $(X,\preceq)$ be an ordered set,
$E\subseteq X$ and $a,b\in E$. A subset $\{e_1,e_2,...,e_k\}$ of $E$
is
called $\prec\succ$-chain between $a$ and $b$ in $E$ if\\
\indent\hspace{0.5mm} (i) $k\geq 2$,\\
\indent\hspace{0.5mm} (ii) $e_1=a~{\rm and}~e_k=b$,\\
\indent\hspace{0.5mm} (iii) $e_1\prec\succ e_2\prec\succ\cdots\prec\succ e_{k-1}\prec\succ e_k$.\\

Let C$(a,b,\prec\succ,E)$ denotes the class of all
$\prec\succ$-chains between $a$ and $b$ in $E$. In particular for
$E=X$, we write C$(a,b,\prec\succ)$ instead of C$(a,b,\prec\succ,X)$.\\

Now, we state and prove some results for uniqueness of coincidence
point, point of coincidence and common fixed point corresponding to
earlier results. For the sake of naturally, firstly we
prove results corresponding to Theorem 3.5 and thereafter for Theorem 3.3.\\

\noindent{\bf Theorem 4.3.} {\it In addition to the hypotheses of
Theorem 3.5,
suppose that the following condition holds:\\
$(u_0)$ C$(fx,fy,\prec\succ,gX)$ is nonempty, for each
$x,y\in X$.\\
Then $f$ and $g$ have a unique point of coincidence.}\\
{\bf Proof.} In view of Theorem 3.5, the set of the coincidence
points (and hence points of coincidence) of $f$ and $g$ is nonempty.
Let $\overline{x}$ and $\overline{y}$ be two points of coincidence
of $f$ and $g$, then $\exists~x,y\in X$ such that
$$\overline{x}=g(x)=f(x)\;{\rm and}\; \overline{y}=g(y)=f(y).\eqno(19)$$
Now, we show that
$$\overline{x}=\overline{y}.\eqno(20)$$
As $f(x),f(y)\in f(X)\subseteq g(X)$, by $(u_0)$, there exists a
$\prec\succ$-chain $\{gz_1,gz_2,...,gz_k\}$ between $f(x)$ and
$f(y)$ in $g(X)$, where $z_1,z_2,...,z_k\in X$. Owing to (19),
without loss of generality, we can choose $z_1=x$ and $z_k=y$. Thus
we have
$$g(z_1)\prec\succ g(z_2)\prec\succ\cdots\prec\succ g(z_{k-1})\prec\succ g(z_k).\eqno(21)$$
Define the constant sequences $z_n^1=x$ and $z_n^k=y$, then using
(19), we have $g(z^1_{n+1})=f(z^1_n)=\overline{x}$ and
$g(z^k_{n+1})=f(z^k_n)=\overline{y}\;\forall~ n\in \mathbb{N}_0.$
Put $z_0^2=z_2,~z_0^3=z_3,...,~ z_0^{k-1}=z_{k-1}$. Since
$f(X)\subseteq g(X)$, on the lines similar to that of Theorem 3.3,
we can define sequences $\{z_n^2\},~ \{z_n^3\},...,~ \{z_n^{k-1}\}$
in $X$ such that $g(z^2_{n+1})=f(z^2_n),~
g(z^3_{n+1})=f(z^3_n),...,~
 g(z^{k-1}_{n+1})=f(z^{k-1}_n)\;\forall~ n\in \mathbb{N}_0.$ Hence,
we have
$$g(z^i_{n+1})=f(z^i_n)\;\;\forall~ n\in \mathbb{N}_0\;{\rm and~for~each}\;i\;(1\leq i\leq k).\eqno(22)$$
Now, we claim that
$$g(z_n^1)\prec\succ g(z_n^2)\prec\succ\cdots\prec\succ g(z_n^{k-1})\prec\succ g(z_n^k)\;\;\forall~ n\in \mathbb{N}_0.\eqno(23)$$
We prove this fact by the method of mathematical induction. Owing to
(21), (23) holds for $n=0.$ Suppose that (23) holds for $n=r>0$,
$i.e.,$
$$g(z_r^1)\prec\succ g(z_r^2)\prec\succ g(z_r^3)\prec\succ\cdots\prec\succ g(z_r^{k-1})\prec\succ g(z_r^k).$$
On using $g$-comparability of $f$, we obtain
$$f(z_r^1)\prec\succ f(z_r^2)\prec\succ f(z_r^3)\prec\succ\cdots\prec\succ f(z_r^{k-1})\prec\succ f(z_r^k),$$
which on using (22), gives rise
$$g(z_{r+1}^1)\prec\succ g(z_{r+1}^2)\prec\succ g(z_{r+1}^3)\prec\succ\cdots\prec\succ g(z_{r+1}^{k-1})\prec\succ g(z_{r+1}^k).$$
It follows that (23) holds for $n=r+1$. Thus, by induction, (23) holds for all $n \in \mathbb{N}_0$.\\
Now for all $n \in \mathbb{N}_0$, define
$${\begin{cases}t_n^1:=d(gz_n^1,gz_n^2)\cr
t_n^2:=d(gz_n^2,gz_n^3)\cr
\vdots \\
$$\hspace{0.0in}t_n^{k-2}:=d(gz_n^{k-2},gz_n^{k-1})\cr
t_n^{k-1}:=d(gz_n^{k-1},gz_n^k).\cr\end{cases}}$$ On using (22),
(23) and assumption (vi), it can be easily shown that
$$t_{n+1}^i\leq \alpha t_{n}^i\;\;\forall~ n\in \mathbb{N}_0\;{\rm and~for~each}\;i\;(1\leq i\leq k-1).$$ By induction, for each
$i\;(1\leq i\leq k-1)$, we get
$$t_{n+1}^i\leq \alpha t_{n}^i\leq \alpha^2 t_{n-1}^i\leq\cdots\leq \alpha^{n+1} t_{0}^i$$
so that
$$ t_{n+1}^i\leq \alpha^{n+1} t_{0}^i.$$ Taking the limit as
$n\rightarrow \infty$ on both the sides of above inequality, we
obtain
$$\lim\limits_{n\to\infty}t_{n}^i=0\;{\rm ~for~each}\;i\;(1\leq i\leq k-1).\eqno(24)$$
On using triangular inequality and (24), we obtain
$$d(\overline{x},\overline{y})\leq t_n^1+t_n^2+\cdots+t_n^{k-1}
\to 0\;\; as \;\; n\to\infty$$
$\Longrightarrow \indent\hspace{4cm}\overline{x}=\overline{y}.$\\
Hence (20) is proved. \\

\noindent{\bf Theorem 4.4.} {\it In addition to the hypotheses of
Theorem
4.3, suppose that the following condition holds:\\
$(u_1)$ one of $f$ and $g$ is one-one.\\
Then $f$ and $g$ have a unique coincidence point.}\\
{\bf Proof.} Let $x$ and $y$ be two coincidence points of $f$ and
$g$, then in view of Theorem 4.3, we have
$$g(x)=f(x)=f(y)=g(y).$$
As $f$ or $g$ is one-one, we have
$$x=y.$$

\noindent{\bf Theorem 4.5.} {\it In addition to the hypotheses of
Theorem 4.3,
suppose that the following condition holds:\\
$(u_2)$ $(f,g)$ is weakly compatible pair.\\
Then $f$ and $g$ have a unique common fixed point.}\\
{\bf Proof.} Let $x$ be a coincidence point of $f$ and $g$. Write
$g(x)=f(x)=\overline{x}$, then in view of Lemma 2.9 and $(u_2)$,
$\overline{x}$ is also a coincidence point of $f$ and $g$. It
follows from Theorem 4.3 with $y=\overline{x}$ that
$g(x)=g(\overline{x}),$ $i.e.$, $\overline{x}=g(\overline{x})$,
which yields that
$$\overline{x}=g(\overline{x})=f(\overline{x}).$$
Hence, $\overline{x}$ is a common fixed point of $f$ and $g$. To
prove uniqueness, assume that $x^*$ is another common fixed point of
$f$ and $g$. Then again from Theorem 4.3, we have
$$x^*=g(x^*)=g(\overline{x})=\overline{x}.$$
This completes proof.\\

\noindent{\bf Theorem 4.6.} {\it In addition to the hypotheses of
Theorem 3.3,
suppose that the condition $(u_0)$ (of Theorem 4.3) holds, then $f$ and $g$ have a unique common fixed point.}\\
{\bf Proof.} We know that in a metric space, every compatible pair
is weakly compatible so that $(u_2)$ trivially holds. Hence,
proceeding on the lines of the proof of Theorems 4.3
and 4.5 our result follows.\\

\noindent{\bf Corollary 4.7.} {\it Theorem 4.3 (resp. Theorem 4.6)
remains true if we replace the condition $(u_0)$ by one of the
following conditions (besides retaining rest of the hypotheses):\\
$(u_0^1)$:  $(fX,\preceq)$ is totally
ordered,\\
$(u_0^2)$: $(X,\preceq)$ is $(f,g)$-directed.}\\
{\bf Proof.} Suppose that $(u_0^1)$ holds, then for each pair
$x,y\in X$, we have $$f(x)\prec\succ f(y),$$ which implies that
$\{fx,fy\}$ is a $\prec\succ$-chain between $f(x)$ and $f(y)$ in
$g(X)$. It follows that C$(fx,fy,\prec\succ,gX)$ is nonempty, for
each $x,y\in X$, $i.e.$, $(u_0)$ holds and hence Theorem 4.3 (resp.
Theorem 4.6) is applicable.\\
Next, assume that $(u_0^2)$ holds, then for each pair $x,y\in X$,
$\exists~z\in X$ such that
$$f(x)\prec\succ g(z)\prec\succ f(y),$$
which implies that $\{fx,gz,fy\}$ is a $\prec\succ$-chain between
$f(x)$ and $f(y)$ in $g(X)$. It follows that
C$(fx,fy,\prec\succ,gX)$ is nonempty, for
each $x,y\in X$, $i.e.$, $(u_0)$ holds and hence Theorem 4.3 (resp. Theorem 4.6) is applicable.\\
\bigskip
\section{A Related Fixed Point Result}
On setting $g = I$, the identity mapping on $X$, in Theorem 3.3
(together with Theorems 4.6), we get the following
fixed point result.\\

\noindent{\bf Corollary 5.1.} {\it Let $(X,d,\preceq)$ be an ordered metric space and $f$ a self-mapping on $X$. Suppose that the following conditions hold:\\
\indent\hspace{0.5mm}$(a)$ $(X,d)$ is complete,\\
\indent\hspace{0.5mm}$(b)$ $f$ is comparable,\\
\indent\hspace{0.5mm}$(c)$ either $f$ is continuous or $(X,d,\preceq)$ has {\it TCC} property,\\
\indent\hspace{0.5mm}$(d)$ there exists $x_{0}\in X$ such that
$x_{0}\prec\succ
f(x_{0})$,\\
\indent\hspace{0.5mm}$(e)$  there exists $\alpha\in [0,1)$ such that \\
\indent\hspace{2.5cm}$d(fx,fy)\leq\alpha d(x,y)\;\;\forall ~x,y\in X$ with $x\prec\succ y$.\\
Then $f$ has a fixed point. Moreover, if we add the following one\\
\indent\hspace{0.5mm}$(f)$ C$(fx,fy,\prec\succ)$ is nonempty, for
each
$x,y\in X$.\\
Then we obtain uniqueness of fixed point.}\\

Corollary 5.1 sharpens Theorems 1.1-1.6 in the following
considerations:\\

\noindent $\bullet$ In Theorem 1.1 and Theorem 1.2, the monotonicity
of $f$ can be replaced by comparability of $f$, which is
relatively weaker.\\

$\bullet$ In Theorem 1.6, the property on $(X,d,\preceq)$ embodied
in assumption $(c)$ is very restrictive and can alternately be
replaced by {\it TCC} property on $(X,d,\preceq)$, which is
relatively weaker.\\
$\bullet$ All the hypotheses of Theorems 1.1-1.5 without assumption
$(f)$ guarantee the existence of fixed point and the presence of
assumption $(f)$ ensures the uniqueness of fixed point. Also
assumption $(f)$ of Corollary 5.1 is relatively weaker than each of
assumptions $(f)$
(of Theorems 1.1-1.6).\\

Finally, we furnish an example which demonstrates that the notion of
comparable mapping is an improvement over monotonicity of the map.\\

\noindent{\bf Example 5.2}. Let $X=[-\frac{1}{3},\frac{1}{3}]$.
Then $(X,d,\preceq)$ is an ordered complete metric space under the
usual metric and the natural partial order. Define $f:X\rightarrow
X$ by $f(x)=x^2,$ then $f$ is comparable but not monotone. Also, for
$x,y\in X$ with $x\preceq y$, we have
$$d(fx,fy)=|x^2-y^2|=|x+y||x-y|\leq \frac{2}{3}d(x,y).$$
$i.e.$ $f$ satisfies the contractivity condition $(e)$ of Corollary
5.1. Thus, all the conditions mentioned in Corollary 5.1 are
satisfied. Notice that $f$ has a unique fixed point in $X$ (namely:
$x=0$).\\

\noindent{\bf Acknowledgment:} Both the authors are grateful to an
anonymous learned referee for his valuable suggestions and fruitful
comments.

\end{document}

%% file: prembu.tex
\newcommand{\vs}{\vspace}
\newcommand{\bc}{\begin{center}}
\newcommand{\ec}{\end{center}}